\documentclass[11pt]{article}

\usepackage{amssymb}
\usepackage{amscd}
\usepackage{amsmath}
\usepackage{graphics}

\newtheorem{theorem}{Theorem}
\newtheorem{proposition}[theorem]{Proposition}
\newtheorem{corollary}[theorem]{Corollary}
\newtheorem{lemma}[theorem]{Lemma}

\newenvironment{proof}{{\bf Proof }}{\par}

\newtheorem{exam}[theorem]{Example}

\newtheorem{rem}[theorem]{Remark}

\newenvironment{remark}{\begin{rem}\rm}{\end{rem}}

\numberwithin{theorem}{section}

\newcommand{\pp}{\mathbb{P}}
\newcommand{\qq}{\mathbb{Q}}
\newcommand{\cc}{\mathbb{C}}

\newcommand{\zz}{\mathbb{Z}}

\newcommand{\git}{/\!\!/}

\newcommand{\M}{\mathcal{M}}
\newcommand{\N}{\mathcal{N}}
\newcommand{\K}{\mathcal{K}}

\begin{document}
\pagestyle{plain}

\title{The stringy E-function of the moduli space of
rank 2 bundles over a Riemann surface of genus 3}
\date{}
\author{Young-Hoon Kiem\\
Department of Mathematics\\
Stanford University\\
Stanford, CA 94305\\
kiem@math.stanford.edu\footnote{Address from September 2002:
Department of Mathematics, Seoul National University, San 56-1,
Seoul, 151-747, South Korea}}

\maketitle

\begin{abstract}
We compute the stringy E-function (or the motivic
integral) of the moduli space of rank 2 bundles
over a Riemann surface of genus 3. In doing so, we answer
a question of Batyrev about the stringy E-functions
of the GIT quotients of linear representations.
\end{abstract}
\footnotetext{AMS Classification numbers:  14F, 14J.}

\section{Statement of the main result}

The stringy E-function is an invariant for singular varieties, due
to Kontsevich, Batyrev, Denef and Loeser, which retains useful
information about the singularities (See \cite{Bat, DenLoe1,
DenLoe2, Craw, Loo}).

Let $X$ be a variety with at worst \emph{log-terminal}
singularities, i.e. \begin{itemize}
\item  $X$ is $\qq$-Gorenstein
\item for a resolution of singularities $\rho:Y\to X$ such that the exceptional locus of $\rho$ is a divisor $D$ whose irreducible components $D_1, \cdots, D_r$ are smooth divisors with only normal crossings, we have
$$K_Y=\rho^*K_X+\sum_{i=1}^ra_iD_i$$
with $a_i>-1$ for all $i$, where $D_i$ runs over all irreducible components of $D$. The divisor $K_Y-\rho^*K_X$ is called the \emph{discrepancy
divisor}.\end{itemize}

For each subset $J\subset I=\{1,2,\cdots,r\}$, define
$D_J=\cap_{j\in J}D_j$, $D_{\emptyset}=Y$ and
$D_J^{0}=D_J\setminus \cup_{j\in I\setminus J}D_j$. Then the
stringy E-function of $X$ is defined by
\begin{equation}\label{EstDef}
E_{st}(X;u,v)=\sum_{J\subset I}E(D_J^0;u,v)\prod_{j\in J}\frac{uv-1}{(uv)^{a_j+1}-1}\end{equation} where $$E(Z)=\sum_{p,q}\sum_{k\ge 0} (-1)^kh^{p,q}(H^k_c(Z;\cc))u^pv^q$$
 is the Hodge-Deligne polynomial for a variety $Z$.

The ``change of variable formula" (Theorem 6.27 in \cite{Bat}, Lemma 3.3 in \cite{DenLoe1}) implies that the function $E_{st}$ is independent of the choice of a resolution.  In particular, if $\rho$ is a crepant resolution (i.e. $\rho^*K_X=K_Y$) then $E_{st}(X;u,v)=E(Y;u,v)$.

A projective $\qq$-Gorenstein algebraic variety of dimension $d$ with at worst log-terminal singularities has the Poincar\'e duality
\begin{equation}\label{Poin}
E_{st}(X;u,v)=(uv)^d\, E_{st}(X;u^{-1},v^{-1})\end{equation}
with $E_{st}(X;0,0)=1$. (Theorem 3.7 in \cite{Bat}.)

In this paper, we compute
the stringy E-function of the moduli space $\N$ of rank 2 bundles of even degree
over a Riemann surface of \emph{genus $3$} with fixed determinant.\footnote{See
\cite{Newstead, Seshadri} for general facts about the moduli space.}
Our main result is the following
\begin{theorem}\label{theorem}
$$\begin{array}{ll}E_{st}(\N)&=
\frac{(1-u^2v)^3(1-uv^2)^3-(uv)^4(1-u)^3(1-v)^3}{(1-uv)(1-(uv)^2)}-
\frac{(uv)^2}{2}\big(
\frac{(1-u)^3(1-v)^3}{1-uv}-\frac{(1+u)^3(1+v)^3}{1+uv}
\big)\\
&+2^6(uv)^{5}(1+uv+(uv)^2)(1+(uv)^2)\big(\frac{uv-1}{(uv)^5-1}\big)^2
\end{array}
$$
The stringy Euler number is
$$e_{st}(\N)=\lim_{u,v\to 1}E_{st}(\N)=31\frac{9}{25}.$$
\end{theorem}

The deepest singularities in the moduli space are
the geometric invariant theory (GIT) quotient $sl(2)^3\git SL(2)$
where the action is the diagonal adjoint action. This is a hypersurface
singularity and that makes the genus $3$ case special.
Batyrev asked  (Question 5.5 in \cite{Bat}) the following\\
{\bf Question (Batyrev)}: Let $X$ be a GIT quotient of $\cc^n$ modulo a linear action of $G\subset SL(n)$. Is it true that $E_{st}(X;u,v)$ is a polynomial?

He showed that this is true when $G$ is abelian or finite. A corollary of our
computation is that the answer is NO in general.
\begin{corollary}\label{corollary}
$$\begin{array}{ll}E_{st}(\cc^9\git SL(2))&=
E([\cc^9\git SL(2)]^s)+\frac{(uv)^3(1+uv+(uv)^2)}{1+uv}\\
&+(uv)^{5}(1+uv+(uv)^2)(1+(uv)^2)\big(\frac{uv-1}{(uv)^5-1}\big)^2
\end{array}$$
where $[\cc^9\git SL(2)]^s$ denotes the smooth part of $\cc^9\git SL(2)$.
\end{corollary}
Since $E([\cc^9\git SL(2)]^s)$ is a polynomial, we deduce that
the stringy E-function of $\cc^9\git SL(2)$ is \emph{not} a polynomial.

When the genus of the Riemann surface is 2, the moduli space
is isomorphic to $\pp^3$ and thus the E-function is $1+uv+(uv)^2+(uv)^3$.
When the genus is greater than 3, the deepest
singularities are no longer hypersurface singularities and it doesn't seem
possible to find the discrepancy divisor by explicit computation
as in this paper.

In \S 2, we study the singularities of the moduli space $\N$.
In \S\S 3,4,5, we work out the blow-ups to get a desingularization of $\N$.
We compute the discrepancy divisor in \S 6 and we prove Theorem \ref{theorem}
and Corollary \ref{corollary} in \S 7. We conclude this paper with
a formula for the stringy E-function of the moduli space
$\M$ of rank 2 bundles of even degree, without fixing determinant,
over a Riemann surface of genus 3.

\underline{Acknowledgements}: I would like to thank Tamas Hausel,
Michael Thaddeus, Frances Kirwan and Jun Li for useful
discussions. Also, I am grateful to the referee for comments and
suggestions which led to improvement in exposition.


\section{The moduli space}
The moduli space $\N$ of rank 2 semistable bundles of degree $0$
with trivial determinant over a Riemann surface of genus $g=3$ is
a singular projective variety of complex dimension $6$. The
singularities are Gorenstein by Theorem A of \cite{DN} and
log-terminal as we will see in \S 6. We refer to \cite{Newstead,
Seshadri, KirM} for general results on the moduli space.

The singular locus in $\N$ is the Kummer variety $\K$, which corresponds
to those rank 2 bundles $L\oplus L^{-1}$ for some line bundle $L$ of degree $0$.
The involution $L\to L^{-1}$ gives us a $\zz_2$ action on the Jacobian $Jac_0$
and the Kummer variety $\K$ is identified  with $Jac_0/\zz_2$.
There are $2^{2g}$ fixed points $\zz_2^{2g}=\{[L\oplus L^{-1}]\,
:\, L\cong L^{-1}\}$. Thus we have a stratification
$$\N=\N^s\sqcup (\K -\zz_2^{2g})\sqcup \zz_2^{2g}.$$

The moduli space $\N$ is constructed as the GIT quotient of a smooth
quasi-projective variety $\frak{R}$, which is a subset of the space of
holomorphic maps from the Riemann surface to the Grassmannian $Gr(2,p)$
of $2$-dimensional quotients of $\cc^p$ where $p$ is a large even number,
by the action of $G=SL(p)$. By deformation theory, the slice at a point
$h\in \frak{R}$, which represents $L\oplus L^{-1}$ where $L\cong L^{-1}$, is
$$H^1(End_0(L\oplus L^{-1}))\cong H^1(\mathcal{O})\otimes sl(2)$$
where the subscript $0$ denotes the trace-free part. According to Luna's slice
theorem, there is a neighborhood of the point $[L\oplus L^{-1}]$ with
$L\cong L^{-1}$, isomorphic to $H^1(\mathcal{O})\otimes sl(2)\git SL(2)$
since the stabilizer of the point is $SL(2)$ (\cite{KirM} (3.3)). Because $\dim H^1(\mathcal{O})=g$,
the deepest singularities are just
$$sl(2)^g\git SL(2)=\mathrm{Spec}\, \cc[z_1,\cdots,z_{3g}]^{SL(2)}.$$
By the classical invariant theory (see \cite{Weyl} or more
precisely \cite{Hueb} 5.1), there is an explicit description of
the generators and relations of the invariant subring
$\cc[z_1,\cdots,z_{3g}]^{SL(2)}$. The special feature of the case
$g=3$ is that the quotient $X:=sl(2)^g\git SL(2)$ is a
hypersurface: For each
$(\overrightarrow{u}_1,\overrightarrow{u}_2,
\overrightarrow{u}_3)\in sl(2)^3$, let
$x_1=\overrightarrow{u}_1\cdot \overrightarrow{u}_1,
x_2=\overrightarrow{u}_2\cdot \overrightarrow{u}_2,
x_3=\overrightarrow{u}_3\cdot \overrightarrow{u}_3,
x_4=\overrightarrow{u}_1\cdot \overrightarrow{u}_2,
x_5=\overrightarrow{u}_1\cdot \overrightarrow{u}_3,
x_6=\overrightarrow{u}_2\cdot \overrightarrow{u}_3,
x_7=det(\overrightarrow{u}_1,\overrightarrow{u}_2,\overrightarrow{u}_3)$.
Then $\cc^9\git SL(2)$ is the hypersurface of
$\cc^7=\mathrm{Spec}\,\cc[x_1,x_2,\cdots,x_6,x_7]$ given by the
equation
$$
f(x_1,\cdots,x_7)=x_1x_2x_3+2x_4x_5x_6-x_1x_6^2-x_2x_5^2-x_3x_4^2-x_7^2.
$$
The locus of $\K$ in this neighborhood, as a set, is given by
\begin{equation}\label{Kl}
x_7=0,\ \ x_1x_2-x_4^2=0,\ \ x_1x_3-x_5^2=0,\ \ x_1x_6-x_4x_5=0
\end{equation}
because a point in $\K_X:=\K\cap X$ can be represented by the
$\cc^*$-fixed points
$(\overrightarrow{u}_1,\overrightarrow{u}_2,\overrightarrow{u}_3)$ with
$\overrightarrow{u}_i\in\{diag(a,-a): a\in \cc\}$, $i=1,2,3$.

Next, we consider the middle stratum $\K-\zz_2^{2g}$. Once again,
if we consider a point $h\in \frak{R}$ representing $L\oplus
L^{-1}$ with $L\ncong L^{-1}$, the slice to the orbit is
isomorphic to
\begin{equation}\label{Knormal}
H^1(End_0(L\oplus L^{-1}))\cong H^1(\mathcal{O})\oplus
H^1(L^2)\oplus H^1(L^{-2}).\end{equation}
The stabilizer $\cc^*$ acts with weights $0,2,-2$ respectively on the components.
Hence, there is a neighborhood of the point $[L\oplus L^{-1}]\in \K-\zz_2^{2g}$
in $\N$, isomorphic to
$$H^1(\mathcal{O})\bigoplus \big(H^1(L^2)\oplus
H^1(L^{-2})\git \cc^*  \big).$$
Notice that $H^1(\mathcal{O})$ is the tangent space to $\K$ and hence
$$
H^1(L^2)\oplus H^1(L^{-2}) \git \cc^* =\cc^{2g-2}\git\cc^*$$ is
the normal cone. The GIT quotient of the projectivization $\pp
\cc^{2g-2}$ by the induced $\cc^*$ action is
$\pp^{g-2}\times\pp^{g-2}$ and the normal cone
$\cc^{2g-2}\git\cc^*$ is  obtained by collapsing the zero section
of the line bundle
$\mathcal{O}_{\pp^{g-2}\times\pp^{g-2}}(-1,-1)$.

\section{First blow-up}

We will desingularize the moduli space $\N$ by blowing up three times.
In this section, we describe the first blow-up.

Let $\N_1$ be the blow-up of $\N$ along the deepest strata $\zz_2^{2g}$
and $D'_1$ be the exceptional divisor.
Since the deepest singularities are all $X:=\cc^9\git SL(2)$, we consider
only one of them. The GIT quotient $X$ is the hypersurface of $\cc^7$ with
the equation
$$
f(x_1,\cdots,x_7)=x_1x_2x_3+2x_4x_5x_6-x_1x_6^2-x_2x_5^2-x_3x_4^2-x_7^2.
$$
We blow up at the origin and denote the exceptional divisor also by $D'_1$.
In terms of a local chart, the blow-up map is
\begin{equation}\label{chart1}
(y_1,\cdots,y_7)\to (y_1, y_1y_2,\cdots, y_1y_7).
\end{equation}
We have $f(x_1,\cdots,x_7)=y_1^2g_1(y_1,\cdots, y_7)$ where
\begin{equation}\label{eqn1}
g_1(y_1,\cdots,y_7)=y_1(y_2y_3+2y_4y_5y_6-y_6^2-y_2y_5^2-y_3y_4^2)-y_7^2.
\end{equation}
Hence, the blow-up $X_1$ is the hypersurface given by $g_1$
and the exceptional divisor $D'_1$ is the subset $y_1=0, y_7=0$ in the local
chart. Let $\tilde{\K}_X$
be the proper transform of $\K_X$.

The singular set of $X_1$ in this chart is, by solving $\nabla g_1=0$,
the union of
\begin{equation}\label{delta1}
y_1=0,\ \ \ y_7=0,\ \ \ y_2y_3+2y_4y_5y_6-y_6^2-y_2y_5^2-y_3y_4^2=0
\end{equation}
and
\begin{equation}\label{K1loc}
y_7=0,\ \  y_2-y_4^2=0,\ \ y_3-y_5^2=0,\ \ y_6-y_4y_5=0.
\end{equation}
Notice that the second component of the singular set is just
the proper transform $\tilde{\K}_X$ in view of (\ref{Kl}).

Now we switch to other charts. Since $x_1,x_2,x_3$ are symmetric,
we consider, for instance,
\begin{equation}\label{chart5}
(y_1,\cdots,y_7)\to (y_5y_1, y_5y_2,y_5y_3,y_5y_4,y_5,y_5y_6, y_5y_7).
\end{equation}
In this chart, $X_1$ is given by the equation
\begin{equation}\label{eqn5}
g_5(y_1,\cdots,y_7)=y_5(y_1y_2y_3+2y_4y_6-y_1y_6^2-y_2-y_3y_4^2)-y_7^2.
\end{equation}
and $D'_1$ by $y_5=0, y_7=0$. The singular locus in this chart is the union
of \begin{equation}\label{delta2}
y_5=0,\ \ \ y_7=0,\ \ \ y_1y_2y_3+2y_4y_6-y_1y_6^2-y_2-y_3y_4^2=0
\end{equation}
and
\begin{equation}\label{K5loc}
y_7=0,\ \  y_1y_2-y_4^2=0,\ \ y_1y_3-1=0,\ \ y_1y_6-y_4=0.
\end{equation}
Again the second component is $\tilde{\K}_X$ by comparing with
(\ref{Kl}).\footnote{There is one more chart $(y_1,\cdots,y_7)\to
(y_1y_7, \cdots, y_6y_7, y_7)$ but it doesn't intersect with the
exceptional divisor.}

From the local descriptions
(\ref{delta1}), (\ref{delta2}), we see that the first component
of the singular set is the subvariety
$$
\Delta_X=\{(y_1:\cdots:y_7)\,|\, y_7=0,
y_1y_2y_3+2y_4y_5y_6-y_1y_6^2-y_2y_5^2-y_3y_4^2=0\}
$$
of the projective space $\pp^6$.
\begin{lemma}
$\Delta_X\cong \pp^2\times\pp^2/\zz_2.$
\end{lemma}
\begin{proof}
Define a morphism $\pp^2\times \pp^2\to \pp^6$ by
$$\big( (x:y:z), (p:q:r) \big)\to (2xp:2yq:2zr:xq+yp:xr+zp:yr+zq:0).$$
In fact, this came from the identity
$$\begin{array}{ll}
(xt+ys+zu)(pt+qs+ru)&=2xp\frac{t^2}{2}+2yq\frac{s^2}{2}+2zr\frac{u^2}{2}\\
&+(xq+yp)ts+(xr+zp)tu+(yr+zq)su.
\end{array}$$
Since $\cc[t,s,u]$ is a UFD, the morphism is a 2:1 map whose image is
precisely $\Delta_X$ as one can easily check.
\end{proof}
\vspace{.7cm} The singular locus of $X_1$ is thus $\Delta_X\cup
\tilde{\K}_X$. By direct computation, the singular locus $\pp^2$
of $\Delta_X$ is the intersection $\Delta_X\cap \tilde{\K}_X$
which is the exceptional divisor of the proper transform
$\tilde{\K}_X\to \K_X$. For instance, in terms of the local chart
of (\ref{chart1}), the singular locus of (\ref{delta1}) is given
by the equations of (\ref{delta1}) and (\ref{K1loc}).

We denote by $\Delta$ the disjoint union of $2^6$ $\Delta_X$'s in
the exceptional divisor in $\N_1$ which has $2^6$ components. Then
$\N_1$ is smooth away from $\Delta \cup \tilde{\K}$.

\section{Second blow-up}
In this section, we consider the second blow-up. Namely, we blow up
$\N_1$ along the proper transform $\tilde{\K}$ of $\K$. This is
particularly important because it is the partial desingularization
of $\N$, defined in \cite{KirP}.

Let $\N_2$ be the blow-up of $\N_1$ along $\tilde{\K}$. Let
$D'_2$ be the exceptional divisor and $\tilde{D}'_1$ be the
proper transform of $D'_1$ which has $2^{2g}$ connected components.
We will describe $\N_2$ as the partial desingularization of $\N$.
For more details on partial desingularization,
we refer to \cite{KirM} and \cite{KirP}.

Let $H$ be a reductive subgroup of $G=SL(p)$ and define $Z^{ss}_H$
as the set of semistable points in $\frak{R}$ fixed by $H$. Let
$\frak{R}_1$ be the blow-up of $\frak{R}^{ss}$ along the smooth subvariety
$GZ^{ss}_{SL(2)}$. Then by Lemma 3.11 in \cite{KirP},
the GIT quotient $\frak{R}^{ss}_1\git G$ is
the first blow-up $\N_1$. The $\cc^*$-fixed point set in $\frak{R}_1^{ss}$
is the proper transform $\tilde{Z}^{ss}_{\cc^*}$ of $Z_{\cc^*}^{ss}$ and
the quotient of $G\tilde{Z}^{ss}_{\cc^*}$ by $G$ is $\tilde{\K}$.
If we denote by $\frak{R}_2$ the blow-up of $\frak{R}^{ss}_1$ along the
smooth subvariety $G\tilde{Z}^{ss}_{\cc^*}=G\times
_{N^{\cc^*}}\tilde{Z}^{ss}_{\cc^*}$ where $N^{\cc^*}$ is the normalizer
of $\cc^*$, the GIT quotient $\frak{R}_2^{ss}\git G$ is our second blow-up
$\N_2$ again by Lemma 3.11 \cite{KirP}. This is Kirwan's partial
desingularization of $\N$ (See \S 3 \cite{KirM}), which is an orbifold.

By applying the algorithm for Betti numbers described in
\cite{KirP}, the Poincar\'e series $P(\N_2)=\sum_{k\ge 0}t^k \dim
H^k(\N_2)$ can be computed as follows. By \cite{KirL}, the
equivariant Poincar\'e series $P^G(\frak{R}^{ss})=\sum_{k\ge 0}t^k
\dim H^k_G(\frak{R}^{ss})$ is given by the gauge theoretic
computation of Atiyah and Bott in \S 11 of \cite{AB} and we get
$$P^G(\frak{R}^{ss})=\frac{(1+t^3)^6-t^8(1+t)^6}{(1-t^2)(1-t^4)}.$$
In order to get $\frak{R}^{ss}_1$ we blow up $\frak{R}^{ss}$ along
$GZ_{SL(2)}^{ss}$ and delete the unstable strata. So we get
$$P^G(\frak{R}_1^{ss})=P^G(\frak{R}^{ss})+2^6\big(\frac{t^2+t^4+\cdots+t^{16}}{1-t^4}
-\frac{t^{10}(1+t^2+t^4)}{1-t^2}\big).$$ Now $\frak{R}_2^{ss}$ is
obtained by blowing up $\frak{R}_1^{ss}$ along
$G\tilde{Z}^{ss}_{\cc^*}$ and deleting the unstable strata. Thus
we have
\begin{equation}\begin{array}{ll}
P^G(\frak{R}_2^{ss})=P^G(\frak{R}_1^{ss})
&+(t^2+t^4+t^6)\big(\frac12 \frac{(1+t)^6}{1-t^2}+\frac12
\frac{(1-t)^6}{1+t^2} +2^6\frac{t^2+t^4}{1-t^4}\big)\\
&-\frac{t^4(1+t^2)}{1-t^2}\big((1+t)^6+2^6(t^2+t^4)\big).
\end{array}
\end{equation}
Because the stabilizers of the $G$ action on $\frak{R}^{ss}_2$ are
all finite, we have $$H^*_G(\frak{R}_2^{ss})\cong
H^*(\frak{R}^{ss}_2/G)=H^*(\N_2)$$ and hence we deduce that
\begin{equation}\label{N2Pr}
\begin{array}{ll}
P(\N_2)&=\frac{(1+t^3)^6-t^8(1+t)^6}{(1-t^2)(1-t^4)}\\
&+2^6\big(\frac{t^2+t^4+\cdots+t^{16}}{1-t^4}
-\frac{t^{10}(1+t^2+t^4)}{1-t^2}\big)\\
&+(t^2+t^4+t^6)\big(\frac12 \frac{(1+t)^6}{1-t^2}+\frac12
\frac{(1-t)^6}{1+t^2} +2^6\frac{t^2+t^4}{1-t^4}\big)\\
&-\frac{t^4(1+t^2)}{1-t^2}\big((1+t)^6+2^6(t^2+t^4)\big).
\end{array}
\end{equation}
See \cite{KirM} for the Betti number computation of  the partial desingularization of $\M$, the moduli
space without fixing determinant.

Furthermore, we can refine the above computation to get the
Hodge-Deligne polynomial for $\N_2$ since the observation in \S 14
\cite{KirT} tells us that the morphisms involved in the above
Betti number computation are strictly compatible with the mixed
Hodge structures. By the gauge theoretic computation of \cite{AB},
the Hodge-Deligne series for the equivariant cohomology
$H^*_G(\frak{R}^{ss})$ is\footnote{The recent article \cite{EK} by
Earl and Kirwan contains detailed arguments for the Hodge number
computation of the equivariant cohomology.},
$$\frac{(1-u^2v)^3(1-uv^2)^3-(uv)^4(1-u)^3(1-v)^3}{(1-uv)(1-(uv)^2)}.$$
Blowing up along $GZ^{ss}_{SL(2)}$
and deleting unstable part amounts to adding
$$2^6\big(\frac{uv+(uv)^2+\cdots+(uv)^8}{1-(uv)^2}
-\frac{(uv)^5(1+uv+(uv)^2)}{1-uv}\big)$$
and blowing up along $G\tilde{Z}^{ss}_{\cc^*}$ and deleting unstable
points amounts to adding
$$\begin{array}{ll}
&(uv+(uv)^2+(uv)^3)\big(\frac12 \frac{(1-u)^3(1-v)^3}{1-uv}+\frac12
\frac{(1+u)^3(1+v)^3}{1+uv} +2^6\frac{uv+(uv)^2}{1-(uv)^2}\big)\\
&-\frac{(uv)^2(1+uv)}{1-uv}\big((1-u)^3(1-v)^3+2^6(uv+(uv)^2)\big).
\end{array}$$
Therefore, we get
\begin{equation}\label{N2P}
\begin{array}{ll}
E(\N_2)&=\frac{(1-u^2v)^3(1-uv^2)^3-(uv)^4(1-u)^3(1-v)^3}{(1-uv)(1-(uv)^2)}\\
&+2^6\big(\frac{uv+(uv)^2+\cdots+(uv)^8}{1-(uv)^2}
-\frac{(uv)^5(1+uv+(uv)^2)}{1-uv}\big)\\
&+(uv+(uv)^2+(uv)^3)\big(\frac12 \frac{(1-u)^3(1-v)^3}{1-uv}+\frac12
\frac{(1+u)^3(1+v)^3}{1+uv} +2^6\frac{uv+(uv)^2}{1-(uv)^2}\big)\\
&-\frac{(uv)^2(1+uv)}{1-uv}\big((1-u)^3(1-v)^3+2^6(uv+(uv)^2)\big).
\end{array}
\end{equation}
Notice that (\ref{N2P}) reduces to (\ref{N2Pr}) if we put
$u=v=-t$.

In this context, $D'_1$ is the disjoint union of $2^6$ copies of
$\pp(sl(2)^3)\git SL(2)$ and $\tilde{D}'_1$ is its
partial desingularization. The algorithm in \cite{KirP} gives us
\begin{equation}\label{D1PD}
E(\tilde{D}'_1)=2^6(1+uv+(uv)^2)(1+uv+(uv)^2+(uv)^3).
\end{equation}

The normal bundle to $G\tilde{Z}^{ss}_{\cc^*}$ has rank $2g-2=4$
as we saw in (\ref{Knormal}). As $G\tilde{Z}^{ss}_{\cc^*}\cong
G\times_{N^{\cc^*}}\tilde{Z}^{ss}_{\cc^*}$ from \cite{KirP} and
the normal bundle can be written similarly, the quotient of the
normal bundle by $G$ is the quotient of its restriction to
$\tilde{Z}^{ss}_{\cc^*}$ by the action of $N^{\cc^*}$. If we first
take the quotient by the identity component $N^{\cc^*}_0$ of
$N^{\cc^*}$, we get a $\cc^{4}\git\cc^*$-bundle over
$\tilde{Jac}$, the blow-up of $Jac$ along $\zz_2^{6}$, since
$\tilde{Z}_{\cc^*}\git N^{\cc^*}_0\cong \tilde{Jac}$. Hence there
is a neighborhood of $\tilde{\K}$ in $\N_1$, which is isomorphic
to the $\zz_2$-quotient of the $\cc^{4}\git\cc^*$-bundle over
$\tilde{Jac}$ because $\pi_0(N^{\cc^*})=\zz_2$. As we mentioned at
the end of \S2, the normal cone $\cc^4\git\cc^*$ is obtained by
collapsing the zero section of the line bundle
$\mathcal{O}_{\pp^{1}\times\pp^{1}}(-1,-1)$ and thus the
exceptional divisor $D_2'$ is the $\zz_2$ quotient of the
$\pp^{1}\times\pp^{1}$ bundle over $\tilde{Jac}$. Hence, the
E-polynomial of $D'_2$ is
\begin{equation}\label{D2HN}\begin{array}{ll}
E(D_2')&=[\big((1-u)^3(1-v)^3+2^6(uv+(uv)^2)\big)(1+uv)^2]^{\zz_2}\\
&=\big(\frac12 (1-u)^3(1-v)^3+\frac12 (1+u)^3(1+v)^3 +2^6
(uv+(uv)^2)\big)(1+uv+(uv)^2)\\
&+\big(\frac12 (1-u)^3(1-v)^3-\frac12 (1+u)^3(1+v)^3 \big)(uv)
\end{array}
\end{equation}
where $[\cdot]^{\zz_2}$ denotes the $\zz_2$-invariant part.
The intersection of the two divisors $D'_2$ and $\tilde{D}'_1$ has
$2^6$ components, each of which is isomorphic to a bundle over $\pp^2$
with fiber $\pp^2=\pp^1\times_{\zz_2}\pp^1$.

Now, we can compute the E-function of the smooth part
$\N^s=\N-\K=\N_2-D_2'-\tilde{D}_1'$. The E-polynomial of $\N_2$,
$\tilde{D}_1'$ are (\ref{N2P}), (\ref{D1PD}) respectively.
The E-polynomial of $D_2'-\tilde{D}'_1$ is
$$\begin{array}{ll}
&\big(\frac12 (1-u)^3(1-v)^3+\frac12 (1+u)^3(1+v)^3 -2^6\big)(1+uv+(uv)^2)\\
&+\big(\frac12 (1-u)^3(1-v)^3-\frac12 (1+u)^3(1+v)^3 \big)(uv).
\end{array}
$$
by subtracting
$E(D_2'\cap\tilde{D}'_1)=2^6(1+uv+(uv)^2)^2$ from (\ref{D2HN}).
Therefore, the E-polynomial of $\N^s$ is
\begin{equation}\label{Nsmooth}\begin{array}{ll}
E(\N^s)&=E(\N_2)-E(\tilde{D}_1')-E(D_2'-\tilde{D}'_1)\\
&=\frac{(1-u^2v)^3(1-uv^2)^3-(uv)^4(1-u)^3(1-v)^3}{(1-uv)(1-(uv)^2)}\\
&-\frac12(\frac{(1-u)^3(1-v)^3}{1-uv}+\frac{(1+u)^3(1+v)^3}{1+uv}).
\end{array}\end{equation}

To end this section, we consider the singular locus of $\N_2$. At
a point in $D'_2\setminus\tilde{D}'_1$, $\N_2$ looks like a line
bundle over $\pp^1\times\pp^1$ times $\cc^3$ and hence smooth. The
singular locus thus lies in $\tilde{D}_1'$ and so we restrict our
concern to $X_1$, the blow-up of $X=sl(2)^3\git SL(2)$. We know
from the previous section that $X_1$ is smooth at points in
$D'_1\setminus\Delta_X$. Hence, the singular locus of $\N_2$ lies
over $\Delta$. We claim that the proper transform $\tilde{\Delta}$
of $\Delta$ is precisely the singular locus in $\N_2$. To verify
our claim, we return to the local chart description.

In terms of the local chart (\ref{chart1}), $X_1$ is given by the equation
(\ref{eqn1}) and $\tilde{\K}_X$ is by (\ref{K1loc}). We introduce
new coordinates
$$w_1=y_1,\ w_2=y_2-y_4^2,\ w_3=y_3-y_5^2,\ w_4=y_4,\ w_5=y_5,\
w_6=y_6-y_4y_5, w_7=y_7.$$
Then the equation of $X_1$ is $w_1(w_2w_3-w_6^2)-w_7^2$ and
$\tilde{\K}_X$ is given by $w_2=w_3=w_6=w_7=0$. The blow-up along
$\tilde{\K}_X$ can be now described locally as
\begin{equation}\label{bl12}
(t_1,\cdots, t_7)\to (t_1, t_2, t_2t_3, t_4, t_5, t_2t_6, t_2t_7).
\end{equation}
Since $w_1(w_2w_3-w_6^2)-w_7^2=t_2^2\big(
t_1(t_3-t_6^2)-t_7^2\big)$ in this chart, $X_2$ is given by the equation
\begin{equation}\label{g12}
g_{12}(t_1,\cdots, t_7)=t_1(t_3-t_6^2)-t_7^2.\end{equation}
The singular locus is, from $\nabla g_{12}=0$,
\begin{equation}\label{g12delta}
t_1=0,\ \ t_7=0,\ \ t_3-t_6^2=0\end{equation}
which is the proper transform of $\Delta_X$ in view of the fact that
$\Delta_X$ is $w_1=w_7=0,\ w_2w_3-w_6^2=0$ from (\ref{delta1}).
Similarly, one can use other charts for the second blow-up to check
that the proper transform $\tilde{\Delta}_X$ of $\Delta_X$ is the singular
locus over the local chart (\ref{chart1}).

In the local chart (\ref{chart5}), $X_1$ is given by (\ref{eqn5})
and $\tilde{\K}_X$ is by (\ref{K5loc}) while $\Delta_X$ is given
by (\ref{delta2}). Since we are interested in a neighborhood of
$\Delta_X\cap \tilde{\K}_X$ where $y_1\ne 0$, we may assume that
$y_1\ne 0$. We introduce new coordinates
$$\begin{array}{ll}
&w_1=y_1,\ w_2=y_2-y_4^2/y_1,\ w_3=y_3-1/y_1,\ w_4=y_4,\\
 &w_5=y_5,\ \
w_6=y_6-y_4/y_1, \ \ w_7=y_7.\end{array}$$
In terms of $w$-coordinates, $X_1$ is just $w_1w_5(w_2w_3-w_6^2)-w_7^2$
and $\tilde{\K}_X$ is $w_2=w_3=w_6=w_7=0$. The blow-up map along $\tilde{\K}_X$ can be written locally as (\ref{bl12}) for instance. One can check
again that the singular locus of $\N_2$  over the local chart (\ref{chart5})
is precisely $\tilde{\Delta}_X$.

By a similar computation for each local chart for $X_1$, we deduce
that the singular locus of $\N_2$ is $\tilde{\Delta}$ as claimed. Observe
from the above that $\tilde{\Delta}$ and $\tilde{D}'_1$ are smooth.

\section{Third blow-up}
To obtain a desingularization $\tilde{\N}$ of $\N$, we blow up $\N_2$ along
$\tilde{\Delta}$. Let $D_3$ be the exceptional divisor of this third blow-up and $D_1, D_2$ denote the proper transforms of $\tilde{D}'_1, D_2'$ respectively.

\includegraphics{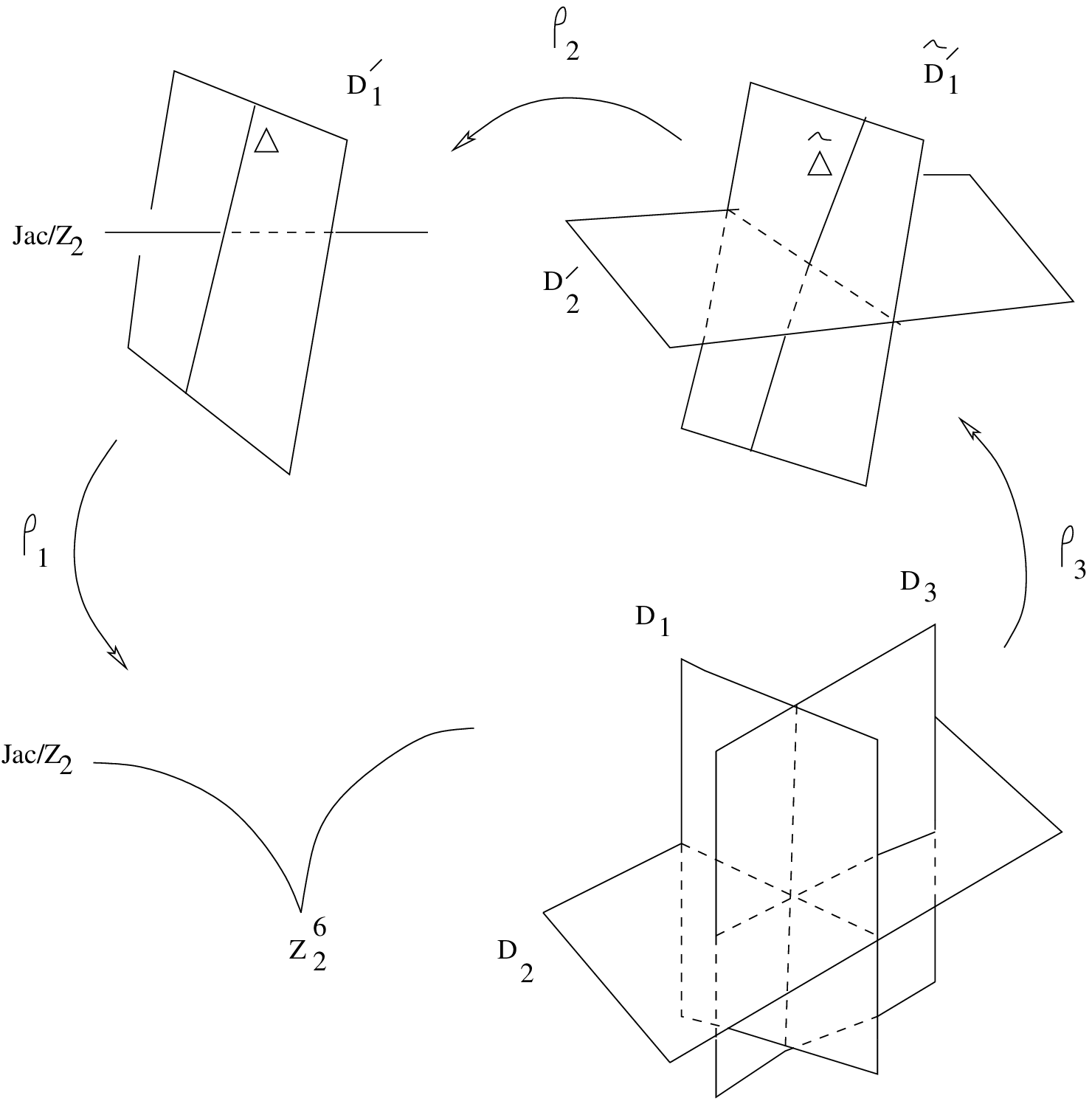}

In terms of the $t$-coordinates (\ref{bl12}) of $\N_2$, one can readily
deduce from (\ref{g12}) and (\ref{g12delta}) that
the singularity along $\tilde{\Delta}$ is just the $(xy=z^2)$-singularity in $\cc^3$ and
by blowing up along $\tilde{\Delta}$ we get a smooth variety. As one can check, the same is true for each local chart of $\N_2$.
Hence, $\tilde{\N}$ is smooth.

One can also explicitly check in terms of local coordinates that
the divisors $D_1, D_2, D_3$ are smooth divisors with only normal crossings.
For instance, consider the $t$-coordinates (\ref{bl12}) for $\N_2$ again. Before blowing up, we introduce new coordinates
$r_1=t_1, r_2=t_2, r_3=t_3-t_6^2, r_4=t_4, r_5=t_5, r_6=t_6, r_7=t_7$.
Then $\N_2$ is given by $r_1r_3-r_7^2=0$ and the blow-up center
is $r_1=r_3=r_7=0$. If we consider the local description of the third blow-up
for instance
$$(\alpha_1,\cdots,\alpha_6,\alpha_7)\to (\alpha_3\alpha_1,\alpha_2,\alpha_3,\alpha_4, \alpha_5, \alpha_6,\alpha_3\alpha_7),$$
$\tilde{\N}$ is $\alpha_1-\alpha_7^2=0$, $D_1$ is $\alpha_1=\alpha_7=0$,
$D_2$ is $\alpha_2=0=\alpha_1-\alpha_7^2$ and $D_3$ is $\alpha_3=0
=\alpha_1-\alpha_7^2$. By repeating a simliar computation for each chart, we
see that the divisors have only normal crossings.

The desingularization process we described can be schematically
summarized in the above picture.

\section{Canonical divisors}
The purpose of this section is to prove the following.

\begin{proposition}\label{discr}
 If $\rho:\tilde{\N}\to \N$ is the desingularization described above, then $K_{\tilde{\N}}=\rho^*K_{\N}+4D_1+D_2+4D_3.$\end{proposition}

\vspace{.7cm}

We consider a differential $$s=\frac{dx_1\wedge dx_2\wedge dx_3\wedge dx_5\wedge dx_6\wedge dx_7}{\partial f/\partial x_4}$$
on $X=sl(2)^3\git SL(2)$. On the smooth part of $X$, $s$ is not vanishing and thus the divisor of $s$ is zero. (See (1.7) \cite{Rei}.)
In terms of local coordinates, the first blow-up map $\rho_1$ is given by $(y_1,\cdots,y_6,y_7)\to (y_1,y_1y_2,\cdots,y_1y_6,y_1y_7)$
and we have a rational differential on $X_1$
$$s=y_1^4\frac{dy_1\wedge dy_2\wedge dy_3\wedge dy_5\wedge dy_6\wedge dy_7}{\partial g_1/\partial y_4}$$
where $f(x_1,\cdots,x_7)=y_1^2g_1(y_1,\cdots,y_7)$. Hence,
$K_{\N_1}=\rho_1^*K_{\N}+4D'_1$.

Now we switch to the $w$-coordinates
$w_1=y_1, w_2=y_2-y_4^2, w_3=y_3-y_5^2, w_4=y_4, w_5=y_5, w_6=y_6-y_4y_5$. Then $g_1=w_1(w_2w_3-w_6^2)-w_7^2$. The second blow-up,
 in terms of local coordinates, is
$(t_1,\cdots,t_7)\to (t_1,t_2,t_2t_3,t_4,t_5,t_2t_6,t_2t_7)$ and we get
a rational differential on $X_2$
$$\begin{array}{ll}
s&=w_1^4\frac{dw_1\wedge dw_2\wedge dw_3\wedge dw_5\wedge dw_6\wedge dw_7}{\partial g_1/\partial w_3}\\
&=t_1^4t_2 \frac{dt_1\wedge dt_2\wedge dt_3\wedge dt_5\wedge dt_6\wedge dt_7}{\partial g_{12}/\partial t_3}\end{array}$$
where $g_1(y_1,\cdots,y_7)=t_2^2g_{12}(t_1,\cdots,t_7)$. Hence,
$$K_{\N_2}=\rho_2^*\rho_1^*K_{\N}+4\tilde{D}'_1+D'_2.$$

We next use the $r$-coordinates $r_1=t_1, r_2=t_2, r_3=t_3-t_6^2,
r_4=t_4, r_5=t_5, r_6=t_6, r_7=t_7$. Then $g_{12}=r_1r_3-r_7^2$.
Finally, we blow up along $r_1=r_3=r_7=0$. In terms of local
coordinates, the blow-up is $(\alpha_1,\cdots,\alpha_7) \to
(\alpha_3\alpha_1,\alpha_2,\alpha_3,\alpha_4, \alpha_5,
\alpha_6,\alpha_3\alpha_7)$ and
$g_{12}=\alpha_3^2(\alpha_1-\alpha_7^2)$. The equation for
$\tilde{\N}$ in the $\alpha$-coordinates is thus $g_{123}=
\alpha_1-\alpha_7^2$ and we have a rational differential on
$\tilde{X}$
$$s=\alpha_1^4\alpha_2\alpha_3^4\frac{d\alpha_2\wedge d\alpha_3\wedge \cdots\wedge d\alpha_7}{\partial g_{123}/\partial \alpha_1}.$$
By a similar computation for each chart, we deduce that
$$K_{\tilde{\N}}=\rho^*K_{\N}+4D_1+D_2+4D_3.$$

\section{The stringy E-function}
We can now compute the stringy E-function of the moduli space $\N$.

The E-function of the smooth part is from \S 4
$$\begin{array}{ll}
E(\N^s)&=E(\N_2)-E(\tilde{D}_1')-E(D_2'-\tilde{D}'_1)\\
&=\frac{(1-u^2v)^3(1-uv^2)^3-(uv)^4(1-u)^3(1-v)^3}{(1-uv)(1-(uv)^2)}\\
&-\frac12(\frac{(1-u)^3(1-v)^3}{1-uv}+\frac{(1+u)^3(1+v)^3}{1+uv}).
\end{array}$$

Next, $D_1^0=D_1-(D_2\cup D_3)$
is $\tilde{D}'_1-\tilde{\Delta}-\tilde{D}'_1\cap D'_2$. Since each component of
$\tilde{D}'_1\cap D'_2$ is a $\pp^2$-bundle over $\pp^2$ and each component of
$\tilde{\Delta}- D'_2$ is $\pp^2\times_{\zz_2}\pp^2-\pp^2$,
$E(D_1^0)$ is (\ref{D1PD}) minus $E(\pp^2\times\pp^2)$ and
 $2^6[(1+uv+(uv)^2)^2]^{\zz_2}-
2^6(1+uv+(uv)^2)$. Hence,
$$E(D_1^0)\frac{uv-1}{(uv)^5-1}=2^6((uv)^5-(uv)^2)\frac{uv-1}{(uv)^5-1}.$$

Since $D_2^0=D_2-(D_1\cup D_3)=D_2'-\tilde{D}_1'$, the E-function
of $D_2^0$ is (\ref{D2HN}) minus
the E-function of $D'_2\cap \tilde{D}'_1$, $2^6(1+uv+(uv)^2)^2$. Hence,
$$\begin{array}{ll}
E(D_2^0)\frac{uv-1}{(uv)^2-1}
&=\big(\frac12 (1-u)^3(1-v)^3+\frac12 (1+u)^3(1+v)^3
-2^6\big)(1+uv+(uv)^2)\frac{uv-1}{(uv)^2-1}\\
&+\big(\frac12 (1-u)^3(1-v)^3-\frac12 (1+u)^3(1+v)^3 \big)(uv)
\frac{uv-1}{(uv)^2-1}.\end{array}$$

As $D_3\cap D_1$ is isomorphic to $\tilde{\Delta}$
and a component of $D_3\cap D_2$ is
a $\pp^1$-bundle over $\tilde{\Delta}\cap D_2'$,
we see that
the E-function of
$D_3^0$ is $2^6$ times the E-function of a $\pp^1$-bundle
over $\tilde{\Delta}$ minus $E(\tilde{\Delta})$ and
$2^6$ times $E((\pp^1-pt)\times
\pp^1\times\pp^2)$. Hence,
$$E(D_3^0)\frac{uv-1}{(uv)^5-1}=2^6((uv)^3+(uv)^4+(uv)^5)\frac{uv-1}{(uv)^5-1}.$$

Notice that $D_{12}^0=D_1\cap D_2-D_3=\tilde{D}'_1\cap D_2'-\tilde{\Delta}$
is the disjoint union of $2^6$ copies of a
$(\pp^2-\pp^1)$-bundle over $\pp^2$. Hence,
$$E(D_{12}^0)\frac{uv-1}{(uv)^5-1}\frac{uv-1}{(uv)^2-1}=
2^6((uv)^2+(uv)^3+(uv)^4)\frac{uv-1}{(uv)^5-1}\frac{uv-1}{(uv)^2-1}.$$

Also, $D_{13}^0=D_1\cap D_3
-D_2$ is $\tilde{\Delta}$ minus $2^6$ $\pp^1$-bundles over $\pp^2$. Hence,
$$E(D_{13}^0)\big( \frac{uv-1}{(uv)^5-1}\big) ^2=
2^6((uv)^2+(uv)^3+(uv)^4)\big( \frac{uv-1}{(uv)^5-1}\big) ^2.$$

Finally, a component of $D_{23}^0=D_2\cap D_3-D_1$ is a $(\pp^1-pt)$-bundle over a $\pp^1$-bundle over $\pp^2$ and a component of $D_{123}^0=D_1\cap D_2\cap D_3$ is a $\pp^1$-bundle
over $\pp^2$. Therefore,
$$E(D_{23}^0)\frac{uv-1}{(uv)^5-1}\frac{uv-1}{(uv)^2-1}=
2^6(uv+(uv)^2+(uv)^3)\frac{uv-1}{(uv)^5-1}$$
and
$$E(D_{123}^0)\big( \frac{uv-1}{(uv)^5-1}\big) ^2\frac{uv-1}{(uv)^2-1}=
2^6(1+uv+(uv)^2)\big( \frac{uv-1}{(uv)^5-1}\big) ^2.$$

Putting together all the pieces above,
we get from the formula (\ref{EstDef}) that
$$\begin{array}{ll}E_{st}(\N)&=
\frac{(1-u^2v)^3(1-uv^2)^3-(uv)^4(1-u)^3(1-v)^3}{(1-uv)(1-(uv)^2)}-
\frac{(uv)^2}{2}\big(
\frac{(1-u)^3(1-v)^3}{1-uv}-\frac{(1+u)^3(1+v)^3}{1+uv}
\big)\\
&+2^6(uv)^{5}(1+uv+(uv)^2)(1+(uv)^2)\big(\frac{uv-1}{(uv)^5-1}\big)^2
\end{array}
$$
This satisfies the Poincar\'e duality
(\ref{Poin}) which serves as a check for our result. Notice that
it is \emph{not} a polynomial.

To prove Corollary \ref{corollary}, let $D_{j,X}$ be the divisors
in $\tilde{X}$ corresponding to $D_j$. Then from above, we have
$$E(D_{1,X}^0)\frac{uv-1}{(uv)^5-1}=((uv)^5-(uv)^2)\frac{uv-1}{(uv)^5-1}$$
$$E(D_{2,X}^0)\frac{uv-1}{(uv)^2-1} =((uv)^3-1)(1+uv+(uv)^2)\frac{uv-1}{(uv)^2-1}$$
$$E(D_{3,X}^0)\frac{uv-1}{(uv)^5-1}=((uv)^3+(uv)^4+(uv)^5)\frac{uv-1}{(uv)^5-1}$$
$$E(D_{12,X}^0)\frac{uv-1}{(uv)^5-1}\frac{uv-1}{(uv)^2-1}=
((uv)^2+(uv)^3+(uv)^4)\frac{uv-1}{(uv)^5-1}\frac{uv-1}{(uv)^2-1}$$
$$E(D_{13,X}^0)\big( \frac{uv-1}{(uv)^5-1}\big) ^2=
((uv)^2+(uv)^3+(uv)^4)\big( \frac{uv-1}{(uv)^5-1}\big) ^2$$
$$E(D_{23,X}^0)\frac{uv-1}{(uv)^5-1}\frac{uv-1}{(uv)^2-1}=
(uv+(uv)^2+(uv)^3)\frac{uv-1}{(uv)^5-1}$$ and
$$E(D_{123,X}^0)\big( \frac{uv-1}{(uv)^5-1}\big) ^2\frac{uv-1}{(uv)^2-1}=
(1+uv+(uv)^2)\big( \frac{uv-1}{(uv)^5-1}\big) ^2.$$ By putting
them together, we get
$$\begin{array}{ll}E_{st}(\cc^9\git SL(2))&=
E([\cc^9\git SL(2)]^s)+\frac{(uv)^3(1+uv+(uv)^2)}{1+uv}\\
&+(uv)^{5}(1+uv+(uv)^2)(1+(uv)^2)\big(\frac{uv-1}{(uv)^5-1}\big)^2
\end{array}$$
where $[\cc^9\git SL(2)]^s$ denotes the smooth part of $\cc^9\git
SL(2)$.

\begin{remark}
If we denote by $\M$ the moduli space of rank 2 semistable bundles of
even degree over a Riemann surface of genus 3 (without fixing determinant),
 the stringy E-function is
\begin{equation}
\label{lasteq}\begin{array}{ll} E_{st}(\M)&=(1-u)^3(1-v)^3\{
\frac{(1-u^2v)^3(1-uv^2)^3-(uv)^4(1-u)^3(1-v)^3}{(1-uv)(1-(uv)^2)}\\
&-\frac{(uv)^2}{2}\big(
\frac{(1-u)^3(1-v)^3}{1-uv}-\frac{(1+u)^3(1+v)^3}{1+uv}
\big)\\
&+(uv)^{5}(1+uv+(uv)^2)(1+(uv)^2)\big(\frac{uv-1}{(uv)^5-1}\big)^2\}.
\end{array}
\end{equation}
We just sketch the computation and leave the details to the
reader. The determinant map $det:\M\to Jac$ is a fibration with
fiber $\N$ and $\M$ has the same singularities as $\N$. So we need
3 blow-ups exactly as in \S\S 3,4,5 and the discrepancy divisor is
given as in Proposition \ref{discr}. It is now easy to modify the
computation to get
$$\begin{array}{ll}
E(\M^s)
&=(1-u)^3(1-v)^3[\frac{(1-u^2v)^3(1-uv^2)^3-(uv)^4(1-u)^3(1-v)^3}{(1-uv)(1-(uv)^2)}\\
&-\frac12(\frac{(1-u)^3(1-v)^3}{1-uv}+\frac{(1+u)^3(1+v)^3}{1+uv})].
\end{array}$$
$$E(D_{1,\M}^0)\frac{uv-1}{(uv)^5-1}=(1-u)^3(1-v)^3((uv)^5-(uv)^2)\frac{uv-1}{(uv)^5-1}$$
$$\begin{array}{ll}
E(D_{2,\M}^0)\frac{uv-1}{(uv)^2-1}
&=(1-u)^3(1-v)^3\\&[\big(\frac12 (1-u)^3(1-v)^3+\frac12
(1+u)^3(1+v)^3
-1\big)(1+uv+(uv)^2)\\
&+\big(\frac12 (1-u)^3(1-v)^3-\frac12 (1+u)^3(1+v)^3 \big)(uv)]
\frac{uv-1}{(uv)^2-1}\end{array}$$
$$E(D_{3,\M}^0)\frac{uv-1}{(uv)^5-1}=(1-u)^3(1-v)^3((uv)^3+(uv)^4+(uv)^5)\frac{uv-1}{(uv)^5-1}$$
$$E(D_{12,\M}^0)\frac{uv-1}{(uv)^5-1}\frac{uv-1}{(uv)^2-1}=(1-u)^3(1-v)^3
((uv)^2+(uv)^3+(uv)^4)\frac{uv-1}{(uv)^5-1}\frac{uv-1}{(uv)^2-1}$$
$$E(D_{13,\M}^0)\big( \frac{uv-1}{(uv)^5-1}\big) ^2=(1-u)^3(1-v)^3
((uv)^2+(uv)^3+(uv)^4)\big( \frac{uv-1}{(uv)^5-1}\big) ^2$$
$$E(D_{23,\M}^0)\frac{uv-1}{(uv)^5-1}\frac{uv-1}{(uv)^2-1}=(1-u)^3(1-v)^3
(uv+(uv)^2+(uv)^3)\frac{uv-1}{(uv)^5-1}$$ and
$$E(D_{123,\M}^0)\big( \frac{uv-1}{(uv)^5-1}\big)
^2\frac{uv-1}{(uv)^2-1}=(1-u)^3(1-v)^3 (1+uv+(uv)^2)\big(
\frac{uv-1}{(uv)^5-1}\big) ^2.$$ Combining these we get
(\ref{lasteq}).

\end{remark}



\begin{thebibliography}{Bat}

\bibitem[AB]{AB}
M.~Atiyah and R.~Bott.
\newblock Yang-Mills equations over Riemann surfaces.
\newblock {\em Phil. R. Soc. Lond. A} 308, 1982.
\newblock Pages 523--615.

\bibitem[Bat]{Bat}
V.~Batyrev.
\newblock Stringy Hodge numbers of varieties
with Gorenstein canonical singularities.
\newblock In {\em Integrable systems and algebraic geometry (Kobe/Kyoto,1997)}, 1998.
\newblock Pages 1--32.

\bibitem[Cra]{Craw}
A.~Craw.
\newblock An introduction to motivic integration.
\newblock Preprint, math.AG/9911179.

\bibitem[DL1]{DenLoe1}
J.~Denef and F.~Loeser.
\newblock Germs of arcs on singular varieties and motivic integration.
\newblock {\em Invent. Math.} 135, no 1, 1999.
\newblock Pages 201--232.

\bibitem[DL2]{DenLoe2}
J.~Denef and F.~Loeser.
\newblock Motivic Igusa zeta functions.
\newblock {\em J. Algebraic Geometry} 7, 1998.
\newblock Pages 505--537.

\bibitem[DN]{DN}
J.-M.~Drezet and M. Narasimhan.
\newblock Groupe de Picard des vari\'et\'es de modules de
fibr\'es semi-stables sur les courbes alg\'ebriques.
\newblock {\em Invent. Math.} 97, 1989.
\newblock Pages 53--94.

\bibitem[EK]{EK}
R.~Earl and F.~Kirwan.
\newblock The Hodge numbers of the moduli spaces of vector bundles
over a Riemann surface.
\newblock Preprint, math.AG/0012260.

\bibitem[Hue]{Hueb}
J.~Huebschmann.
\newblock Poisson Geometry of flat connections for SU(2)-bundles
on surfaces.
\newblock {\em Math. Zeit.} 221, 1996.
\newblock Pages 243--259.

\bibitem[Ki1]{KirM}
F.~Kirwan.
\newblock On the homology of compactifications of moduli spaces of
vector bundles over a Riemann surface.
\newblock {\em Proc. Lon. Math. Soc.} 53, 1986.
\newblock Pages 237--266.

\bibitem[Ki2]{KirL}
F.~Kirwan.
\newblock On spaces of maps from Riemann surfaces to
Grassmannians and applications to the cohomology of moduli spaces of vector bundles.
\newblock {\em Ark. Mat.} 24, 1986.
\newblock Pages 221--275.

\bibitem[Ki3]{KirP}
F.~Kirwan.
\newblock Partial desingularizations of quotients of nonsingular
varieties and their Betti numbers.
\newblock {\em Ann. Math.} 122, 1985.
\newblock Pages 41--85.

\bibitem[Ki4]{KirT}
F.~Kirwan.
\newblock Cohomology of quotients in algebraic and symplectic geometry.
\newblock Mathematical Notes 31.
\newblock Princeton, 1985.

\bibitem[Loo]{Loo}
E.~Looijenga.
\newblock Motivic measures.
\newblock Preprint, math.AG/0006220.

\bibitem[New]{Newstead}
P.~Newstead.
\newblock Introduction to moduli problems and orbit spaces.
\newblock Tata Institute Lecture Note, 1978.

\bibitem[Rei]{Rei}
M.~Reid.
\newblock Young person's guide to canonical singularities.
\newblock In {\em Algebraic Geometry Bowdoin}.
\newblock Vol.46, {\em Proc. Symp. Pure Math.}, AMS, 1987.

\bibitem[Ses]{Seshadri}
C.~Seshadri.
\newblock Fibr\'es vectoriels sur les courbes alg\'ebriques.
\newblock {\em Ast\'erisque} 96, (1982).

\bibitem[Wey]{Weyl}
H.~Weyl.
\newblock The classical groups.
\newblock Princeton University Press, 1946.



\end{thebibliography}
\end{document}